\documentclass[12pt]{amsart}
\usepackage{latexsym, amsmath, amscd, amssymb, amsthm,natbib} 
\newtheorem{te}{Theorem}[section]
 \newtheorem{oz}{Definition}[section]
 
 \newtheorem{lm}{Lema}[section] 
 
\begin{document} 

\title{ Casimir elements and kernel of Weitzenb$\ddot{\rm{\bf{O}}}$k derivation. }
 \author[L. Bedratyuk]{L. Bedratyuk}
\noindent
 \address{ Khmelnitskiy national university \newline 
\noindent
Insituts'ka st., 11\\ 
\noindent
Khmelnitskiy, 29016\\
\noindent
 Ukraine}\email{bedratyuk@ief.tup.km.ua } 
\begin{abstract} Let  $k[X]{:=}\,k[x_0,x_1,\ldots ,x_n]$ be a polynomial algebra over a  field $k$ of characteristic zero. We offer an algorithm for calculation of kernel of Weitzenb$\ddot{\rm{o}}$k derivation ${d(x_i)=x_{i-1}},  \ldots,$ ${ d(x_0)=0}$, ${i=1\ldots n}$  that is based on an analogue of the well known Casimir elements of finite dimensional Lie algebras. By using this algorithm, the kernel  is calculated in the case $n <7$.
 \end{abstract}
 \maketitle
 \section{Introduction}
\noindent
Let  $k[X]{:=}\,k[x_0,x_1,\ldots ,x_n]$ be a polynomial algebra over a  field $k$ of characteristic zero. For arbitrary derivation   $D$ of   $k[X]$ denote by   $k[X]^D$  a kernel of  $D$, that is,
$$ k[X]^{D}:={\rm{Ker} }D:=\{f\in k[X]; D(f)=0 \}. $$
\noindent
Let $d$ be  Weitzenb$\ddot{\rm{o}}$k derivation defined by rule ${d(x_i)=x_{i-1}, \ldots, d(x_0)=0}$, ${i=1\ldots n}.$ The kernel  of the  derivation $d$ was actively studied by various
 authors. The well-known Weitzenb$\ddot{\rm{o}}$k's theorem   implies a finite generation of the derivation $d$. A minimal generating set of algebra $k[X]^{\,d}$ for $n\leq 4$  presented in   $[1]$ and in  $[5]$ for  $n\leq 5.$
The aim of this paper is to calculate a generating set of the kernel  of  Weitzenb$\ddot{\rm{o}}$k derivation in case $n\leq 6.$
We offer the general  description of a kernel of arbitrary polynomialy derivation $D$ by using a construction which is a commutative analogue to a construction of  Casimir  elements of finite dimension Lie algebras. Let us recal that a (generalised) Casimir element of a finite dimensional  Lie algebra $L$ is called an cental element of the universal eneveloping algebra $U(L)$ of the follow form
$$ \sum_{i} u_i\, u_i^* ,$$
\noindent
where  $\{u_i\}$, $\{u_i^*\}\ $ \! are a dual bases of a contragradient  $L\,-$ modules in  $U(L)$ with respect to the adjoint action of the Lie algebra on the algebra $U(L)$. In the case ${\mbox{char}(k){=}0}$ it is well known that every  element of the center of an universal eneveloping algebra $U(L)$  is  a Casimir element. The case  ${\mbox{char}(k){=}p{>0}}$ was studied by present autor in  $[4]$.

\noindent
This paper is organized as follows.
\noindent
In section 2 a conceptions of $D\,-$ module, dual $D\,-$\! modules and Casimir elements are introduced. We prove that arbitrary Casimir elements of derivation  $D$ belongs to the kernel  $k[X]^{\,D}$. 
For any linear derivation $D$ is showed that  a theorem inverse to above  theorem  is true, i.e., any element of kernel $k[X]^{\,D}$ is a Casimir element.

\noindent
In section 3 we study Casimir elements(polynomial) of the  Weitzenb$\ddot{\rm{o}}$k derivation $d$. Since   $d$ is linear derivation we see that the problem of finding of the kernel  $k[X]^{\,d}$ is equivalent to the problem of finding  a realisation of a dual $d\,-$\! modules in  $k[X]$. We get such realisations from any element of the kernel $k[X]^{\,d}$ by using two new derivations $e$ and  $\hat  d$ which arising with the natural embedding of   $d\,-$\! modules into  ${\mathfrak sl_2}\,-$\! modules.

\noindent
In section 4 we introduce a maps $\tau_i:k[X]^{\,d}\longrightarrow k[X]^{\,d}$. To any element  $z$ of kernel derivation $d$ we assign certain family  $\tau_i(z)$ elements of the kernel and studying   properties of this correspondence.

\noindent
In section 5 we present a criterion to verify if a subalgebra of $k[X]^{\,d}$ coincides with the whole algebra $k[X]^{\,d}$. By using this criterion we offer an algorithm for computing the kernel $k[X]^{\,d}$.

\noindent
By using the algoritm in section 6 we compute and present  a list of generating elements of the kernel of of Weitzenb$\ddot{\rm{o}}$k derivation in the case $n <7$. The result for $n=6$ is new.

\noindent
All objects which we use in this paper:   rings,  algebras,  vector spases,  isomorphisms,  polynomials are considered over field $k$ of a characteristic zero.

 \section{Casimir's elements of a derivation.}
 \noindent
 The aim of this section is  to offer a method of  constructing  elements of $k[X]^D$ where $D$ is an arbitrary derivation of $k[X]$.
 \begin{oz} A vector space $V$ is called a $D\:-$ module if ${D(V) \subseteq V}$. \end{oz}
 \noindent
 The derivation $d$ defined by $d(x_n)=x_{n-1}$, $d(x_{n-1})=x_{n-2}$,\ldots, $d(x_0)=0$ is called the Weitzenb$\ddot{\rm{o}}$k derivation. Let $X_m$ be a vector space spanned by the elements $x_0,x_1,\ldots,x_m$. Then, obviously, $X_m$ is $d\,-$ module. We denote by $D_V$ a matrix of derivation $D$ in some ordered fixed basis of the space $V$. For example the matrix of derivation $d$ in $X_m$ is Jordan cell $J_{m+1}(0).$
 \begin{oz}
 {$D\,-$ module }\! $V$ is called dual to $D\,-$ module $V^{\,*}$ if there are dual bases $\{v_i\}$, $\{v_i^*\}$ of $V$ and $\,V^*$ such that $$D\,_{V^*}=(-D\,_V)^T.$$
 \end{oz}
\noindent
 The bases $\{v_i\}$, $\{v_i^*\}$ are also called    dual bases.
From definition 2.2 it follows that the derivation $d$ acts on $X_m^*$ such that  ${d(x_i^*)=-x_{i+1},d(x_m)=0}.$

\noindent
Let $D, \hat D$ be a derivations of $k[X]$.
 \begin{oz}
 $D\,-$ module $V$ and $\hat D\,-$ module $W$ are called isomorphic if there is a linear space isomorphism $\phi:V \longrightarrow W$ such that $D\,\phi=\phi\, \hat D$. In this case we write $V \cong W$.
 \end{oz}
\noindent

 \begin{te}\label{dual} $X_m^{\,*}\cong X_m$  as $d\,-$ modules.
 \end{te}
 \begin{proof}
 \noindent Let $ \varphi:X_m^* \longrightarrow X_m $ be the isomorphism defined by rule $ \varphi(x_i^*)=(-1)^{i+1}x_{m-i+1}.$
 \noindent 
Then
 $$d(\varphi(x_i^*))=d((-1)^{i+1}x_{m-i+1})=(-1)^{i+1}x_{m-i},$$ and  $$\varphi(d(x_i^*))=\varphi(-x_{i+1}^*)=-(-1)^{i+2}x_{m-i}=(-1)^{i+1}x_{m-i}= d(\varphi(x_i^*)).$$ 
\noindent
 Hense $ \varphi  $ is the isomorphism from $d\,-$ module $X_m^*$ to $d\,-$ module $\, X_m$ and bases $\{x_i\}$, $\{(-1)^{i+1}x_i \}$ are dual ones. 
\end{proof}
 \begin{oz}\label{OzCa} 
Suppose  $V=\{v_i\}$, $V^*=\{v_i^*\}$ are dual $D\,-$ modules in $k[X]$. The the polynomial $$\Delta(V,V^*):=\sum_i(v_i \cdot v_i^*).$$ \noindent is called  the  { \bf Casimir element } of a derivation $D$. 
\end{oz}
 \noindent By using theorem \ref{dual} we obtain following Casimir elements of degree 2 for the Weitzenb$\ddot{\rm{o}}$k derivation $d$:
$$ \Delta(X_k,X_k^*)=\sum_{i=1}^k (-1)^{i+1}x_i \cdot x_{k-i+1}.$$

 \noindent
 It is easy to show that $\Delta(X_k,X_k^*) \in k[X]^{\, d}$. In generally the following theorem holds:

 \begin{te}\label{Ca} 
Suppose  $U$ and $U^*$ are two dual $D\,-$ modules in $k[X]$. Then ${\Delta(U,U^*)\in k[X]^D}$.
 \end{te}
 \begin{proof} 
Assume that the bases $\{u_i\}$, $\{u_i^*\}$ are dual ones and ${D_U=\{\lambda_{i\,j}\}}$, ${i,j=1\ldots n}$. Then
 $$
D(u_i)=\sum_{j=1}^{n} \lambda_{i\,j}u_j, D(u_i^*)=\sum_{j=1}^{n}( -\lambda_{j\,i}\,u_j).
$$
 Therefore
 $$
 D(\Delta(U,U^*))=D\Bigl(\sum_{i=1}^{n} u_i u_i^*\Bigr)=\sum_{i=1}^{n} (D(u_i)u_i^*+u_i\,D(u_i^*))=
 $$
 $$ 
=\sum_{i=1}^{n} \Bigl(\sum_{j=1}^{n}(\lambda_{i\,j} u_j)u_i^*+u_iD(u_i^*)\Bigr)=\sum_{i=1}^{n} \Bigl(\sum_{j=1}^{n} \lambda_{j\,i}u_i\,u_j^*+u_i\,D(u_i^*)\Bigr)=
 $$
 $$ 
=\sum_{i=1}^{n}  u_i\Bigl(\sum_{j=1}^{n}\lambda_{j\,i}u_j^*+D(u_i^*)\Bigr)=0. 
$$
 \end{proof}
\noindent
 For any Casimir element one can show that this element does not depend on the choice of dual bases; so that a Casimir element  is well defined.

\noindent
 For linear derivation $D$ of $k[X]$ a theorem inverse to theorem 2.2 is true. Let us assume that 
$$
D(x_i)=\sum_{j=0}^{n} \lambda_{i\,j} x_j.
$$
 \begin{te}\label{main} 
Let $z$ be homogeneous polynomial belonging to $k[X]^D$. Then $z$ is a Casimir element. 
\end{te}
 \begin{proof}
Let us remember that a derivation of the algebra $k[X]$ of a form ${f_0\partial_0+f_1\partial_1+\cdots +f_n\partial_n}$, $f_i \in k[X]$ is called a special derivation. It is well-known that the set of all special derivation $W_n:={\rm Der}(k[X])$ is a Lie algebra with respect to the commutator of a derivations. Taking into account $D=D(x_0) \partial_0 +D(x_1) \partial_1 +\cdots + D(x_{n}) \partial_n$ we get
 $$
[D,\partial_i]=-\sum_{i=0}^{n} \lambda_{j\,i}\partial_{j}.
$$
 \noindent 
It is clear that vector space $X_n$ is $D\,-$\! module. Without loss of generality, we can assume that $\partial_n(z) \neq 0.$
\begin{lm}
 Let $z$ be non-vanishing element from $k[X]^D$. Then the vector space {          }\\ ${Z_D:=\langle \partial_0(z),\cdots ,\partial_n(z) \rangle}$ is $D\,-$ module and $Z_D \cong X_n^*$.
 \end{lm}
 \begin{proof} We need to verify duality of the bases elements now. In fact 
$$ D(\partial_i(z))=[D,\partial_i](z)+\partial_i(D(z))=[D,\partial_i](z)=-\sum_{i=0}^{n} \lambda_{j\,i}\partial_{j}(z). $$
 \end{proof}
 \noindent
 Now since the modules $X_n$ and  $Z_D$ are dual we can write  their Casimir element. By using Euler's theorem about homogeneous polynomials we get 
$$
 \Delta(X_n,Z_D)=x_0\partial_0(z)+x_1\partial_1(z)+\cdots +x_n \partial_n(z)=\mbox{deg}(z)\,z. 
$$
 \noindent
 Thus $z=\frac{{\displaystyle1}}{{\displaystyle \mbox{deg}(z)}} \Delta(X_n,Z_D)$, i.e., $z$ is a Casimir element.
 \end{proof} 
\noindent
We end this section with the following

\noindent
{\bf Conjecture.} {\it For any derivation  $D$  the algebra  $k[X]^{\,D}$  is generated by Casimir elements.}

\section{Casimir elements of Weitzenb$\ddot{\rm{o}}$k derivation.}

\noindent
From theorem  \ref{main} it follow that to know the ring  $k[X]^{\,D}$ we have to know a realisations of {$D\,-$ modules} in $k[X].$  Below we offer a way to construct these realisations for Weitzenb$\ddot{\rm{o}}$k derivation $d$.

\begin{te}

Any  $d\,-$ module $V=\langle v_0,v_1,\ldots v_{n} \rangle $   can be extended to  ${\frak sl}_2\,-$ module  where ${\frak sl}_2\,$  is simple three-dimensional Lie algebra over field $k$.

\end{te}

\begin{proof}
Let us introduce on  $V$ two additional derivations  $ \hat d $ and  $e$ as follows:
$$
\begin{array}{l}
\hat d(v_i):=(i+1)(n-i)v_{i+1},\\
e(v_i):=(n-2\,i)v_i.
\end{array}
$$
\noindent
 By straightforward calculation for any $i$ we get:
$$
\begin{array}{l}
[d,\hat d](v_i)=e(v_i), \\

[d,e](v_i)=-2\,d(v_i),\\

[\hat d,e](v_i)=2\,\hat d(v_i).
\end{array}
$$
\noindent
Those commutator relations coincide with the commutator relations of bases elements of the simple three-dimensional Lie algebra
 ${\frak sl}_2$.  Hence the vector space  $V_n$ together with operators  $d,\hat d, e$  is  ${\frak sl}_2\,-$\! module.
\end{proof}
\noindent 
\begin{oz}
For any polynomial  $z \in k[X]$  a natural number $s$ is called an  {\bf order}  of the polynomial $z$ if the number $s$ is the  smallest natural number such that 
$$
\hat d^s(z) \ne 0, \hat d^{s+1}(z) = 0. 
$$
\end{oz}
\noindent
We denote an order of $z$ by $\mbox{ord}(z)$. For example $\mbox{ord}\,(x_0)=n.$
By using Leibniz's  formula we get 
$$\mbox{ord}(a\cdot b)=\mbox{ord}(a)+\mbox{ord}(b) \mbox{ for all  } a,b \in k[X]. $$
For the  derivation  $\,e\,$ every monomial  $x_0^{\alpha_1}\, x_1^{\alpha_2}\, \cdots x_n^{\alpha_n}$  is an eigenvector with the eigenvalue  \\ $w(x_0^{\alpha_0}\, x_1^{\alpha_1}\, \cdots x_n^{\alpha_n})$  where
$$
w(x_0^{\alpha_0}\, x_1^{\alpha_1}\, \cdots x_n^{\alpha_n})=n\,(\sum_i\, \alpha_i)-2\,(\alpha_1+2\,\alpha_2+\cdots +n\,\alpha_n).
$$
A homogeneous  polynomial is called {\bf isobaric} if all its monomials have equal eigenvalue.
\begin{oz}
An eigenvalue of arbitrary monomial of a homogeneous isobaric  polynomial $z$ is called the   {\bf weight} of the polinomial $z$ and denoted by $\omega (z).$
\end{oz}
\noindent
It is easy to see that $\omega(a\cdot b)=\omega(a)+\omega(b)$ for  a homogeneous isobaric  polynomials $a, b.$ 
\begin{te}
For arbitrary homogeneous isobaric polynomial  ${z\in k[X]^d}$ a vector space  
$$  
 V_m(z){:=}\langle v_0,v_1,\ldots v_{m} \rangle, v_i=\frac{(w(z)-i)!}{i!\,(\omega(z))!}\, \hat d^{\,{i}}(z), v_0:=z, m=0\ldots {s},
$$  
\noindent
is  $d\,-$\! module, moreover $X_m\cong V_m(z) $. 
Here  $\omega(z), s $\!  are weight and order of  $z$.
\end{te}
\begin{proof}
\noindent
Let us  prove two intermediate relations :

\begin{enumerate}
 \item[({\it i})]  $e(\hat d^{\,i}(z))=(\omega(z)-2\,i)\,\hat d^{\,i}(z),$
\item[({\it ii})] $d(\hat d^{\,i}(z))=i\,(\omega(z)-i+1)\,\hat d^{\,i-1}(z). $
\end{enumerate}
\noindent
The relation  $(i)$ is valid for  $i=0$:
$$e(\hat d^{\,0}(z))=e(z)=\omega(z)\, z.$$
\noindent
If this relation  holds for certain $i$, then 
$$
e(\hat d^{\,i+1}(z))=[e,\hat d ](\hat d^{\,i}(z))+\hat d(e(\hat d^{\,i}(z)))=-2\, \hat d^{\,i+1}(z)+\hat d((\omega(z)-2i)\hat d^{\,i+1}(z))=
$$
$$
=(\omega(z)-2(i+1))\,\hat d^{\,{i+1}}(z).
$$
\noindent
The relation  $(ii)$ is true for  $i=0$.  In fact 
$$ d(\hat d^{\,0}(z))=d(z)=0 $$
If this relation  holds for certain  $i$ then by using $(i)$   we get 
$$ 
d(\hat d^{\,{i+1}}(z))=[d, \hat d](\hat d^{\,i}(z)) +\hat d(d(\hat d^{\,i}(z)))=  e(\hat d^{\,i}(z))+\hat d(i(\omega(z){-i}+1)\,\hat d^{\,i-1}(z))=
$$
$$
=(\omega(z)-2\,i)\hat d^{\,i}(z)+i(\omega(z){-}i+1)\hat d^{\,i}(z))= (i+1)(\omega(z)-i)\,\hat d^{\,i}(z). 
$$
\noindent
Hence the relations are valid for any $i$.
\noindent
Consider now a vector space
$$
 V_m(z):=\langle v_0,v_1,\cdots , v_{m} \rangle,
$$ 
 where $ v_i=\alpha_i(z) \, \hat d^{\,i}(z)$ for some undefined  ${ \alpha}_i(z) \in k$. 
For the vector space $V_m(z)$ to be a $d\,-$ module it is enought to have $d(v_i)=v_{i-1}$ for all $i$.
Since
$$
d(v_i)=d(\alpha_i(z) \hat d^{\,i}(z))=\alpha_i(z) \, i\, (\omega(z){-}i)\, \hat d^{\,i-1}(z),
$$
\noindent
we get a following recurrence formula for   $\alpha_i(z)$:
$$ i\,(\omega(z){-}i+1)\,\alpha_i(z)=\alpha_{i-1}(z), \alpha_1(z)=1 ,$$
\noindent
by solving it we obtain 
$$
\alpha_i(z)=\frac{(\omega(z)-i)!}{i! \omega(z)!}.
$$
\end{proof}

\section{Maps $\tau_i$.}

\noindent
For arbitrary homogeneous isobaric polynomial  ${z\in k[X]^d}$ denote by
 $\tau_i(z)$  the Casimir  element: 
$$
\tau_i(z):=\Delta(X_i, V^{\,*}_i(z)), 0\leq i\leq \mbox{min}(\mbox{ord}(z),n), z\in k[X]^{\,d}.
$$
\noindent
Since every polynomial of  $k[X]^{\,d}$ is a sum of a homogeneous isobaric polynomials,  the map  $\tau_i$ can be extended to whole algebra  $k[X]^{\,d}$ by requiring
$ \tau_i(u+v){:=}\tau_i(v)+\tau_i(v).$
\noindent
Besides we clearly   have  $ \tau_i(\lambda z){=}\lambda \tau_i(z), \lambda \in k$. Thus  ${\tau_i:k[X]^{\,d}  \longrightarrow     k[X]^{\,d}}$ is now a linear map. Note that the maps   $\tau_i$ is well defined only on elements $z$ of a kernel such that   $\mbox{ord}(z)\leq i$.

\begin{te}
Any polynomial of   $z \in k[X]^{\,d}$ has the form 
$$
z=\frac{1}{\deg(z)}(\tau_n(c(0))+\tau_{n-1}(c_{1})+\ldots +\tau_0(c(n))),
$$ 
where 
$$
\begin{array}{l}
c(0)=\partial_n(z),\\
c(i)=\partial_{n-i}(z)+\sum\limits_{k=1}^i (-1)^{k+1} c_k(i-k),\\
c_k(i):=\displaystyle\frac{(w(c(i))-k)!}{k!\, w(c(i))!}\,\hat d^k(c(i)), c_0(i):=c(i).
\end{array}
$$
\end{te}
\begin{proof}
We may assume that $z$ is a homogeneous isobaric polynomial. \noindent
We first show that  $d(c(i))=0$ for all $i$.
Since  $[\partial_n$,$d]=0$  we have  $d(c(0))=d(\partial_n(z))=\partial_n(d(z))=0$.
Suppose by induction  $d(c(i))=0$. We have to show that  $d(c(i+1))=0.$ From theorem 3.2 it follows  that  $d(c_k(i))=c_{k-1}(i)$.
By definition we have 
$$
c(i+1)=\partial_{n-(i+1)}(z)+\sum\limits_{k=1}^{i+1} (-1)^{k+1} c_k(i+1-k)=\partial_{n-(i+1)}(z)+c_1(i) +\sum\limits_{k=2}^{i+1} (-1)^{k+1} c_k(i+1-k).
$$
Therefore
$$
d(c(i+1))=d(\partial_{n-(i+1)}(z)) +d(\sum\limits_{k=2}^{i+1} (-1)^{k+1} c_k(i+1-k))+d(c_1(i))=
$$
$$
=-\partial_{n-1}(z)+\sum\limits_{k=2}^{i+1} (-1)^{k+1} c_{k-1}(i+1-k)+c(i)=
$$
$$
=-(\partial_{n-i}(z)+\sum\limits_{k=1}^{i+1} (-1)^{k+1} c_k(i+1-k))+c(i)=-c(i)+c(i)=0.
$$
\noindent
Now since  $\partial_n(z)=c(0)$ and  ${\mbox{ord}(\partial_n(z))=\mbox{ord}(z)+n \geq n}$  then the Casimir element   $\tau_n(c(0))$ exists and 
$$
\tau_n(c(0))=x_n c(0)-x_{n-1}\, c_1(0)+\ldots +(-1)^n x_0 c_n(0). 
$$
\noindent
Furthermore, taking into account  $\deg(z) z=x_n \partial_n(z)+x_{n-1} \partial_{n-1}(z)+\ldots +x_0 \partial_0(z)$  we get :
$$
\deg(z) z-\tau_n(c(0))=x_{n-1} (\partial_n(z)+c_1(0))+x_{n-2} (\partial_{n-2}(z)-c_2(0))+\ldots +x_0 (\partial_0-(-1)^n c_n(0). 
$$
The coefficient of  $x_{n-1}$ is equal   $c(1)$ hence 
$$
\tau_{n-1}(\partial_n(z)+c_1(0))=x_{n-1} c_0(1)-x_{n-2}\, c_1(1)+\ldots +(-1)^n x_0 c_{n-1}(1).
$$
\noindent
Therefore
$$
\deg(z)\,z -(\tau_n(c(0))+\tau_{n-1}(c(1)))=x_{n-2}(\partial_{n-2}(z)-c_2(0)+c_1(1))+\ldots +
$$
$$
+x_i\,(\partial_i(z)-(-1)^i c_i(0)-(-1)^{i-1} c_{i-1}(1) 
+\ldots + x_0\,(\partial_0(z)-(-1)^n c_n(0)-(-1)^{n-1} c_{n-1}(1)=
$$
$$
=x_{n-2}(c(2))+\ldots + x_0\,(\partial_0(z)-(-1)^n c_n(0)-(-1)^{n-1} c_{n-1}(1).
$$
Finally, we obtain
$$
\deg(z)\,z -(\tau_n(c(0))+\tau_{n-1}(c(1))+\ldots+\tau_1(c(n-1))=
$$ 
$$
=x_0\,(\partial_0(z)+c_1(n)-c_2(n-1)+\ldots+(-1)^{n+1} c_n(0)=x_0 \,c(n)=\tau_0(c(n)).
$$
\noindent
Hence
$$
\deg(z) z=\tau_n(c(0))+\tau_{n-1}(c(1))+\ldots +\tau_1(c(n-1))+\tau_0(c(n)),
$$
\noindent
and  we get
$$
z=\frac{1}{\deg(z)}(\tau_n(c(0))+\tau_{n-1}(c_{1})+\ldots +\tau_0(c(n))).
$$

\end{proof}

\begin{te}
Let  $z$ be a  homogeneous isobaric polynomial of  $k[X]^{\,d}$. Then
\begin{enumerate}
  \item[(i)]     ${\rm{ord}}(\tau_i(z))=n+{\rm{ord}}(z){-}2\,i $,
 \item[(ii)]   $  \omega(z)={\rm{ord}}(z)$.

 \end{enumerate}
\end{te}
\begin{proof}
$(i)$  Let us    denote by $s$ the  order of $z$.  We first  show that  $\hat d^{n+s-2\,i}(\tau_i(z){ \neq} 0$  but  \\ $\hat d^{n+s-2\,i+1}(\tau_i(z){=} 0$.
Consider now  two    $\hat d\, -$\! modules: 
$$
\hat X_m:=\langle x_n, \frac{x_{n-1}} {\gamma_{n-1}},  \frac{x_{n-2}} {\gamma_{n-1}\,\gamma_{n-2}},\cdots  \frac{x_m} {\gamma_{n-1}\,\gamma_{n-2}\cdots \gamma_{n-m} }   \rangle,
$$
$$
\hat V_m(z) :=\langle z_s, \frac{z_{s-1}} {\alpha_{s-1}(z)},  \frac{z_{s-2}} {\alpha_{s-1}(z)\,\alpha_{s-2}(z)}\cdots  \frac{z_m} {\alpha_{s-1}(z)\,\alpha_{s-2}(z)\cdots \alpha_{s-m}(z) }   \rangle,
$$
where   $\gamma_i{=}(i+1)\,(n-i)$, $z_i{:=}\alpha_i(z) \hat d^i(z)$, $ z_0{:=}z$, $\alpha_i{=}\displaystyle \frac{(w(z)-i)!}{i! w(z)!}$, $m{=}0..\mbox{min}(n,s)$ and 
$w(z)\,$ is a  weight of the polynomial $z$.
\noindent
Define   a  linear multiplicative map $\psi :X_m\cdot V_m(z) \longrightarrow \hat X_m \cdot \hat V_m(z)$ by the rule 
$$
\psi(x_i)=\frac{x_{n-i}}{\gamma_{n-1}\,\gamma_{n-2} \ldots \gamma_{n-i}}, \psi(z_i)=\frac{z_{n-i}}{\alpha_{n-1}(z)\,\alpha_{n-2}(z) \ldots \alpha_{n-i}(z)}.
$$

\noindent
Since  $\hat x_i=\gamma_i\, x_{i+1}$ and 

$$
\psi(d(x_i))=\psi(x_{i-1})=\frac{x_{n-(i-1)}}{\gamma_{n-1}\,\gamma_{n-2} \ldots \gamma_{n-(i-1)}},
$$
$$
\hat d (\psi(x_i))=\hat d \Bigl(\frac{x_{n-i}}{\gamma_{n-1}\,\gamma_{n-2} \ldots \gamma_{n-i}}\Bigr)=\frac{x_{n-(i-1)}}{\gamma_{n-1}\,\gamma_{n-2} \ldots \gamma_{n-(i-1)}},
$$

\noindent
it follows that  $\psi(d(x_i)){=}\hat d (\psi(x_i))$  thus  the restriction of  $\psi$ to  $X_m$ is a isomorphism from  $d\,-$ module $X_m$  to  $\hat d\,-$ module $\hat X_m$. Similarly     $V_m(z)$  and $\hat V_m(z)$ are isomorphic thus  $\psi$ is a isomorphism of   $d\,-$  module  $X_m\cdot V_m(z)$ to  $\hat d\, -$ module $\hat X_m \cdot \hat V_m(z)$.
\noindent
From Theorem 3.2 it follows that $\hat d\,-$ modules  $\hat X_m$   and   $\hat V_m(z)$   are isomorphic.  For arbitrary homogeneous isobaric polynomials $z\in k[X]^{\,d}$ denote by  $\hat \tau_i(z)$ the  corresponding Casimir element  

$$
\hat \tau_i(z)=\Delta(\hat X_i,\hat V^{\*}_i(z)).
$$
\noindent
It is easy to check  that  $\hat \tau_i(z){=}\psi(\tau_i(z))$. 
\noindent
Let us show that 

$$
\hat d^{\, n+s-2}(\tau_1(z))=\frac{(n+s-2)! (n-1)! (s-1)!}{n \, s}\,  \hat \tau_1(z).
$$

\noindent
Since
$$
\gamma_0 \cdot  \gamma_1 \cdots \gamma_{k-1}=(1\cdot n) (2 \cdot (n-1)) \ldots (k \cdot (n-(k-1))=[k!]^2 { n \choose k},
$$
we have
$$
\hat d^k(x_i)=\gamma_i \gamma_{i+1} \ldots \gamma_{i+k-1}  \cdot x_{k+i} =\frac{\gamma_0 \gamma_{1} \ldots \gamma_{i+k-1} }{\gamma_0 \gamma_{1} \cdots \gamma_{i-1}  } x_{k+i} =\frac{[(i+k)!]^2 \displaystyle {n \choose {k+i}}}{[i!]^2 \displaystyle {n \choose i}} x_{k+i}.
$$
\noindent
In particulary   $\hat d^{\,n-1}(x_0){=}[(n-1)!]^2 n \cdot x_{n-1}$, $\hat d^{\,n-1}(x_1){=}\displaystyle \frac{[(n)!]^2}{n} \cdot x_n.$
Likewise,
$$
\hat d^k(z_i)=\frac{[(i+k)!]^2 \displaystyle {s \choose {k+i}}}{[i!]^2 \displaystyle {s \choose i}} x_{z+i}.
$$

\noindent
Taking into account  $\hat d^{\,n+1}(x_0){=}\hat d^{\,s+1}(z_0){=}0$ and $\tau_1(z)=x_0\,z_1-x_1\,z_0$, we obtain

$$
\hat d^{\, n+s-2}(x_0\,z_1-x_1\,z_0)={{n+s-2}\choose n} \hat d^n(x_0) \cdot \hat d^{\,s-2}(z_1)+{{n+s-2}\choose {n-1}} \hat d^{\,n-1}(x_0) \cdot \hat d^{\,s-1}(z_1))-
$$
\newline
$$
-{{n+s-2}\choose {n-1}} \hat d^{\,n-1} (x_1) \cdot \hat d^{\,s-1}(z_0)-{{n+s-2}\choose {n-2}} \hat d^{\,n-2} (x_1) \cdot \hat d^{\,s}(z_0)=
$$
\newline
$$
=x_{n-1} z_s \left({{n+s-2}\choose {n-1}}\cdot [(n-1)!]^2 n \frac{ [(s-1)!]^2}{s} -{{n+2-2}\choose {n-2}} [(n-1)!]^2 [s!]^2\right)+
$$
\newline
$$
+x_n z_{s-1} \left({{n+s-2}\choose {n}} [n!]^2 [(s-1)!]^2 
-{{n+s-2}\choose {n-1}}\cdot [(s-1)!]^2 s \frac{ [(n-1)!]^2}{n}\right) =
$$

$$
=(n+s-2)!\, s!\, (n-1)! \, \, x_{n-1}\, z_s - (n+s-2)!\, n!\, (s-1)!\,  x_n z_{s-1}=
$$

$$
=(n+s-2)! (n-1)! (s-1)!(s \, x_{n-1} z_s - n\, x_n  z_{s-1})=
$$

$$
=\frac{(n+s-2)! (n-1)! (s-1)!}{n \, s} (\frac{x_{n-1}}{n}\,z_s -x_n \frac{z_{s-1}}{s})
=\frac{(n+s-2)! (n-1)! (s-1)!}{n \, s}\hat \tau_1(z). 
$$

\noindent
In the general case, it can be shown
(the proof is routine) that 

$$
\hat d^{\,n+s-2\,i}(\tau_i(z))=\frac{(n+s-2\,i)! (n-i)! (s-i)!}{n (n-1)...(n-(i-1)) s (s-1)...(s-(i-1))} \hat \tau_i(z) \neq 0,
$$

\noindent
whenever  $ \tau_i(z) \neq 0, $
but  $\hat d^{\,n+s-2\,i+1}(\tau_i(z))=\hat d (\hat d^{\,n+s-2\,i}(\tau_i(z)))=0$  since  $\hat \tau_i(z)$ belongs to the kernel of the derivation  $\hat d$,
therefore the order of non-vanished polinomial  $\tau_i(z)$ is equal  $n+s-2\,i$.

$(ii)$
It is a direct corollary  of the relation $d(\hat d^{\,i}(z))=i\,(\omega(z)-i+1)\,\hat d^{\,i-1}(z). $

\end{proof}

\section{An algorithm of computing of the kernel  $k[X]^{\,d}$ }
\noindent
For any subalgebra  $ T \subseteq k[X]^{\,d}$  we write  $\tau(T)$ for  the subalgebra generated by the  elements    $\tau_i(z), z\in T, i \leq \mbox{ord}(z)$. 
\noindent
The map $\tau$ allow  us to arrange the following  iteration process for calculation of the algebra   $k[X]^{\,d}$ :
 for arbitrary subalgebra $B$ of algebra  $k[X]^{\,d}$ denote by  $\overline B$ a subalgebra  generated by elements  $B \bigcup \tau (B)$.
For every integer  $m \geq 0$ define the following sequence  $B_m$ of subalgebras of  $k[X]^{\,d}$
$$
\left\{
\begin{array}{l}
  B_0=k[x_0] \\
  B_m={\overline{B_{m-1}}}.
\end{array}
\right.
$$

\noindent
We get a increasing chain of the subalgebras $B_i$
$$
B_0 \subseteq B_1 \subseteq B_2 \ldots 
$$
\noindent
Now it is possible to state the following  proposition
\begin{te}
There exists some  $k$ such that  $B_{k}=B_{k+1}=k[X]^{\,d}$.
\end{te}

This theorem is a direct corollary of finite generation of the Weitzenb$\ddot{\rm{o}}$k  derivation and of the follow theorem :
\begin{te}
If  $T$ is a subalgebra of  $k[X]^{\,d}$  with the property that  $x_0 \in T$ and  $\tau (T) \subseteq T$, then  $T=k[X]^{\,d}$.
\end{te}
\begin{proof}
It is enought to show that  $k[X]^{\,d} \subseteq T$. The proof it is by induction on the  degree of a polynomial   $z \in T$. Under the condition of the theorem we have   $x_0 \in T$ . Let us assume that the theorem is true for all  polynomials of degree less or equal to  $s$. Suppose  $\deg(z)=s+1$; then from theorem 4.1 it is follows that  $z$ can be expresed as a sum of polynomials    $\tau_i(z^{\,'}_i)$  where $z^{\,'}_i$  are elements of the degree $s$. Then by the induction hypothesis, we have   $z \in T$ therefore  $k[X]^{\,d} \subseteq T$.
\end{proof}
\noindent
For the  realization of the algorithm it is necessary to be able to calculate   algebra $\tau(B_{i})$ knowing generating set for algebra $B_i$. 
\noindent
Let us give a definition of a irreducible polynomials of algebra  $k[X]^{\,d}$. We will say that   $x_0$ is the only irreducible polynomial of depth unity, and that  $x_0, K_1,\cdots ,K_s$  form a complete set of irreducible polynomials of depth $<m$ if every polynomials of  $B_i$,   $i < m$  is a polynomial in  $x_0,K_1,\ldots ,K_s$, over field $k$. Polynomials which are not irreducible we will call as reducible. 
\noindent
A several cases of reducebility  considered in the following theorem.

\begin{te}
Suppose that  $u,\,v \in k[X]^{\,d}$ are irreducible polynomials of depth $i$; then  
\begin{enumerate}
  \item[(i)]      If  ${\rm ord}(u)=0$   then   $\tau_k(u\,v)$ is reducible of the depth $i$;
  \item[(ii)]    $\tau_i(u\, v)$ is reducible for all  $i \leqslant {\rm min}(n,{\rm ord}(v))$;
 \item[(iii)]  Polynomial  $\tau_i(u_1\, u_2 \ldots u_{i+1}) $, $u_i \in  B_i$  is always reducible . 
\end {enumerate}
\end{te}
\begin{proof}

$(i)$ 
Let  $i \leqslant \mbox{ord}(v)$   and $V_i(v)=\langle v_0v_1\ldots ,v_i \rangle .$ If   $\hat d(u)=0$ then  we have \\ ${V_i({u\, v})=\langle u\,v_0, u\,v_1\ldots ,u\,v_i \rangle }$. Thus 
$$
\tau_i(u\, v)=\Delta(X_i,V^{\,*}_i({u\, v}))=\sum_{k=1}^i (-1)^{k+1} x_k\, u\, v_{i-k+1}=u\, \tau_i(v).
$$
$(ii)$ 
Let us show  that for $i \leqslant {\rm min}(n,{\rm ord}(v))$ and for certain  $c_k^{\,'}\in k[X]^{\,d}$  the relation 

$$
\tau_i(u\, v)=u\, \tau_i(v)+\sum_{k=1}^{i-1}\tau_{i-k}(c_k^{\,'}), 
$$

is true. 
In fact
$$
\tau_i(u\,v)-u\,\tau_i(v)=\sum_{k=0}^i (-1)^{i-k} \alpha_{i-k}(u\,v) x_{i-k}\hat d^{\,k}(u\,v)-u\, \left(\sum_{k=0}^i (-1)^{i-k} \alpha_{i-k}(v) x_{i-k}\hat d^{\,k}(v)\right)=
$$
$$
=x_i\,u v+\sum_{k=1}^i (-1)^{i-k} \alpha_{i-k}(u\,v) x_{i-k}\hat d^{\,k}(u\,v)-u\, \left(x_i v+\sum_{k=1}^i (-1)^{i-k} \alpha_{i-k}(v) x_{i-k}\hat d^{\,k}(v)\right)=
$$
$$
=\sum_{k=1}^i (-1)^{i-k} x_{i-k}(\alpha_{i-k}(u\,v) \hat d^{\,k}(u\,v)-\alpha_{i-k}(v) u\, \hat d^{\,k}(v)
$$

The polynomial of the right hand side belongs to the kernel as difference of two elements of the kernel. To conclude the proof, it remains to apply  theorem 4.1

$(iii)$ If among the polynomials  $u_1,\, u_2, \ldots ,u_{i+1} $ there is a polynomial of the order zero then  $\tau_i(u_1\, u_2 \ldots u_{i+1}) $ is reducible  by  $(i)$. If all of them have non-vanishing orders  then the order  of  $u_1\, u_2 \ldots u_{i}$ is  $\geq i$  and   $\tau_i(u_1\, u_2 \ldots u_{i+1}) $ is reducible by  $(ii)$.  
\end{proof}

\begin{oz}
A polynomial $z \in \tau (B_m)$  is called  {\bf acceptable} for algebra   $B_m$ if $z$ is ireducible and   $z  \not \in B_{m}$.
\end{oz}
\noindent
Let  $B_m{=}k[f_1,\ldots ,f_r,g_1,\ldots ,g_s]$ where $\{g_i\}\,$  are acceptable for  $B_{m-1}$. From  {theorem  4.3} it is follows 

\begin{te}
Polynomial  $ \tau_i(f_1^{\,\alpha_1}\,f_2^{\,\alpha_2}\,\ldots \, f_r^{\,\alpha_r}\,g_1^{\,\beta_1}\,\ldots \, g_s^{\,\beta_s})$  can not be an acceptable polynomial for  $B_m$ if any of follow conditions  holds :
\begin{enumerate}
\item[1.] $\sum\limits_k \alpha_k+\sum\limits_k \beta_k > i .$

 \item[2.] $\sum\limits_k \beta_k=0.$

\item[3.] Some of  $f_i, g_k$  have the order equal to zero but  $\alpha_i \neq 0, \beta_k \neq 0$.

\item[4.]   $ {f_1^{\alpha_1}\,f_2^{\alpha_2}\,f_r^{\alpha_r}\,g_1^{\beta_1}\, g_s^{\beta_s}}$ can be expresed  as product of two polynomials one of them has order greater than  $i$.

\end{enumerate}

\end{te}
\noindent

\begin{oz}
  A triple of  integer numbers 
 $$[{\rm{deg}}(z),{\rm{ord}}(z),\displaystyle \frac{{n\,{\rm{deg}}(z)-\omega(z)}}{2}]$$  is called the signature of a polynomial $z$ and denoted by $[z].$
\end{oz}
\noindent
From theorem 4.2  it is follow that  $[u\,v]=[u]+[v]$

\noindent
Now let us offer a verification algorith  if a polinomial  $z \in k[X]^{\,d}$ belongs to  subalgebra   
 $B=k[f_1,f_2,\ldots ,f_m]$  of algebra   $k[X]^{\,d}$.

\begin{enumerate}

\item Set up  the following system of equation

$$[z]=\alpha_1\,[f_1]+\alpha_2\,[f_2]+\cdots+\alpha_m\,[f_m] $$
\item If this system has no an integer positive solutions then  $z$ doesnt belong to subalgebra  $B$.
\item Suppose   $\{\alpha^{(i)}=(\alpha_1^{(i)},\cdots,\alpha_m^{(i)}),i=1..k\}\, $ \! is the  set of all positive integer solutions of the system. Set up the  new system of equations:
$$z=\beta_1\, f^{\alpha^{(1)}}+\beta_2\, f^{\alpha^{(2)}}+\cdots + \beta_k\, f^{\alpha^{(k)}},$$
 where 
$$
f^{\alpha^{(i)}}=f_1^{\alpha_1^{(i)}}\,f_2^{\alpha_2^{(i)}}\,\ldots f_m^{\alpha_m^{(i)}}.
$$
\item If the system has non-vanishing solutions then   $z\in B$ otherwise     $z \not \in B$ 
\end{enumerate}
\noindent
From above we obtain the following  algorithm for  computing of kernel  $k[X]^{\,d}$.
Let  $\{B\}$ be  a generating set of some  algebra  $B$.
\begin{enumerate}
\item $\{B_1\}=\{x_0\}$.
\item Suppose that algebra    $\{B_i\}=\{f_1,\ldots,f_r\} \bigcup \{g_1,\ldots, g_s\}$ is already calculated  where  $\{g_1,\ldots, g_s\}$  are acceptable polynomials for  $B_i$.
\item  Consider a finite sets of elements of $\tau(B_i)$ which could be an acceptable for $B_i$:
$$
\begin{array}{l}
B_i^{(m)}:=\{\tau_k(f_1^{\alpha_1}\,\cdots \,f_r^{\alpha_r}\,g_1^{\beta_1}\, \cdots g_s^{\beta_s}), \sum\alpha_q+\sum\beta_q=m,  m \leq n, k\leq n \}
\end{array}
$$

\item By using  previuos algoritm we compute a set  $H$  of acceptable polynomials of  $B_i^{(m)}$, $m \leq n$.
\item If  $H=\emptyset $ then  $k[X]^{\,d}=B_{\,i}$ else  ${B_{i+1}}=\{f_1,\ldots,f_r\} \bigcup \{g_1,\ldots, g_s\} \bigcup H$. 
\end{enumerate}

\section{ A calculation $k[X]^{\,d}$ for  $n<7$.}
\noindent
Denote by $t$ the variable $x_0$.
\begin{te}\label{B_2}

$B_1=k[t,\tau_2(t),\tau_4(t),\cdots ,\tau_{2[\frac{n}{2}]}(t)], n>2$

\end{te}
\begin{proof}

It is clear  that an aceptable polynomials for $B_1$ can  only be any one of the following
  polynomials  $\tau_i(t), i\leq n$.  It is easy to check  that for odd $i$ we have $\tau_i(t)=0.$ To prove that the set of polynomials $t$, $\tau_2(t),\tau_4(t),\cdots ,\tau_{2[\frac{n}{2}]}(t)$ is the minimal generating set for subalgebra $B_1$ it is enought to prove  that there are not any linear relations for  the polynomials $t^2$, $\tau_2(t),\tau_4(t),\cdots ,\tau_{2[\frac{n}{2}]}(t).$ The proof follows obviously from the fact that no two polynomials have the same orders.
\end{proof}
\noindent
The following computations were all done with Maple.
\newline

\noindent
{\bf{5.1 \rm {\bf n\,=\,1.}}}
\newline

\noindent
Since    $\tau_1(t)=0$ we get   ${\tau(B_1)=0}$. Therefore, we  have
$B_2=B_1$  thus $k[X]^{\,d}=k[t]$.
\newline

\noindent
{\bf{5.2 \rm {\bf n\,=\,2.}}}
\newline

\noindent
 By using theorem 6.1 we have $B_1=k[t,dv]$ where ${dv:=\tau_2(t)=t\,x_2-2 x_1^2}$. Since ${{\rm ord}({dv})=0}$; then $B_1$ has no an acceptable elements. Therefore  $B_2=B_1$ and $k[X]^{\,d}=k[t,dv]$. 
\newline

\noindent
{\bf{5.3 \rm {\bf n\,= 3.}}}
\newline

\noindent
We have $B_1{=}k[t,dv]$ where ${dv=\tau_2(t)}$. Since ${{\rm ord}({dv})=3+3-2\cdot 2=2}$, we obtain
$$
\begin{array}{l}
B_1^{(1)}=\{\tau_1(dv),\tau_2(dv)\},\\
B_1^{(2)}=\{\tau_3(dv^2)\}.
\end{array}
$$
However  by a straightforward calculation we obtain  $\tau_3(dv^2){=}0$ and  $\tau_2(dv){=}0.$ Denote the  remaining element by  ${tr:= \tau_1(dv)}$.  By using an algorithm of the section 5 we get $tr \not \in B_1$ and thus $B_2=k[t,dv,tr]$.
\noindent
Since ${\rm ord}(tr)=3$ we see that an acceptable elements for  $B_2$ can only  be  the polynomial  ${c=\tau_3(tr)}$.  Note  ${\rm ord}(c)=0$ and  $B_2$ has no an elements of order zero. Hence  $c \not \in B_2$ and $B_3=k[t,dv,tr,c]$. 
Since ${{\rm ord}(c)=0}$, then    $B_3$  has no any acceptable elements. Hense $B_4=B_3$ thus  ${k[X]^{\,d}=k[t,dv,tr,c]}$ where
$$
\begin{array}{l}
dv= - 2\,t\,{x_{2}} + {x_{1}}^{2}, \\
tr= - 3\,t\,{x_{1}}\,{x_{2}} + 3\,t^{2}\,{x_{3}} + {x_{1}}^{3} ,\\
c=- 18\,t\,{x_{1}}\,{x_{2}}\,{x_{3}} + 8\,t\,{x_{2}}^{3} + 9\,{x_{
3}}^{2}\,t^{2} + 6\,{x_{1}}^{3}\,{x_{3}} - 3\,{x_{1}}^{2}\,{x_{2}
}^{2}.
\end{array}
$$
\newline

\noindent
{\bf{5.4 \rm {\bf n\,= 4.}}}
\newline

\noindent
$B_1=k[t,d_1,d_2]$, where
$$
\begin{array}{ll}
 d_1=\tau_2(t) & [d_1]=[2,4,4] ,\\
d_2=\tau_4(t) &  [d_2]=[2,0,12].
\end{array}
$$

\noindent
It is easy to see that for subalgebra  $B_1$ only  polynomial  $\tau_i(d_1)$, ${i=1..4}$ can be an acceptable polynomial. By direct calculation we obtain  $\tau_2(d_1)=t\, d_2$  ³ $\tau_3(d_1)=0$.   Put
$$
\begin{array}{ll}
 tr_1=\tau_1(d_1) & [tr_1]=[3,6,8],\\
 tr_2=\tau_4(d_1) &[tr_2]=[3,0,18].
\end{array}
$$
\noindent
The signatures of  $t^3$, $t\, {d}_1$, $t\, {d}_2$   are equal  $[3,8,0]$, $[3,8,4]$, $[3,4,12]$ therefore  $tr_1$, $tr_2$  are not in  $B_1$, hence   ${B_2=k[t,d_1,d_2,tr_1,tr_2].}$
\noindent
Since  ${\rm ord}(d_2)={\rm ord}(tr_2)=0$ we see that an acceptable elements can  only be the follow elements   
$\tau_i(tr_1)$, ${i=1..4}$. Take into account  $\tau_2(tr_1)=0$, $\tau_2(tr_4)=0$ and  
$$
\begin{array}{l}
 \tau_1(tr_1)=d_2\, t^2+d_1^2, \\
\tau_3(tr_2)=d_1\, d_2-t\cdot tr_2,
\end{array}
$$
\noindent
we have  ${\tau(B_2) \subseteq B_2}$, $B_3=B_2$  and  ${k[X]^{\,d}=k[t,d_1,d_2,tr_1,tr_2]}$  where
$$
\begin{array}{l}
d_1= - 2\,t\,{x_{2}} + {x_{1}}^{2},\\
d_2= - 2\,t\,{x_{4}} + 2\,{x_{1}}\,{x_{3}} - {x_{2}}^{2},\\
tr_1= - 3\,t\,{x_{1}}\,{x_{2}} + 3\,{x_{3}}\,t^{2} + {x_{1}}^{3},\\
tr_2=12\,{x_{2}}\,t\,{x_{4}} - 9\,{x_{3}}^{2}\,t - 6\,{x_{1}}^{2}\,{x
_{4}} + 6\,{x_{1}}\,{x_{2}}\,{x_{3}} - 2\,{x_{2}}^{3}.
\end{array}
$$
\newline

\noindent
{\bf{5.5 \rm {\bf n\,=\,6.}}}
 \newline

\noindent
We have  $B_2=k[t,d_1,d_2]$ where 
$$
\begin{array}{l}
{d_1{:=}\tau_2(t)},{           } {[d_1]=[2,6,2],}\\
{d_2{:=}\tau_4(t)}, {[d_2]=[2,2,4].}
\end{array}
$$
\noindent
The following 10 polynomials can be  acceptable polynomials for  $B_2$ :
$$
\begin{array}{l}
 B_2^{(1)}=\{ \tau_i(d_1),i=1..5;\tau_i(d_2),i=1,2 \},\\
{B_2^{(2)}=\{ \tau_i(d_2^2),i=3,4 \},}\\
{B_2^{(3)}=\{ \tau_5(d_2^3)\}.}
\end{array}
$$
\noindent
By direct calculation we obtain 
$$
\begin{array}{ll}
\tau_2(d_1)=-\frac{6}{5}\,t\,d_2, & \tau_5(d_1)=0\\
\tau(d_2)=5\,\tau_3(d_1) & \tau_2(d_2)=-\frac{5}{4} \tau_4(d_1).
\end{array}
$$
\noindent
Make the denotations,
$$
\begin{array}{ll}
tr_1:=\tau_4(d_1), &     [tr_1]=[3,3,6],\\
tr_2:=\tau_3(d_1), &      [tr_2]=[3,5,5],\\
tr_3:=\tau_1(d_1), &     [tr_3]=[3,9,3],\\
p_1:=\tau_4(d_2^2), & [p_1]=[5,1,12],\\
p_2:=\tau_3(d_2^2), & [p_2]=[5,3,11],\\
si_1:=\tau_5(d_2^3), &[si_1]=[7,1,17].
\end{array}
$$
\noindent
Thus  
$$ B_3=k[t,d_1,d_2,tr_{1-3},p_{1-2},si_1].$$
\noindent
By using an algorithm of section 5 one can show that this polynomial system is the  generating set for $B_3$
\noindent
Now write down  polynomials which can be  acceptable for  $B_2$:
$$
\begin{array}{l}
B_2^{(1)}=\{ \tau_3(tr_1),\tau_5(tr_2),\tau_4(tr_3),\tau_5(tr_3) \}\\
B_2^{(2)}=\{\tau_{4-5}(d_2\,tr_1),\tau_5(d_2\,p_2),\tau_5(tr_1^2)\}\\
B_2^{(3)}=\{\tau_5(d_2^2\,p_1),\tau_5(d_2^2\,tr_1),\tau_5(d_2^2\,si_1),\tau_5(tr_1\,p_1\,si_1),\tau_5(p_1\,p_2\,si_1) \}\\
B_2^{(4)}=\{\emptyset \}
\end{array}
$$
\noindent
By direct calculation we obtain that the  following polynomials are equal to zero : $\tau_2(tr_1)$,\\$\tau_5(d_2\,tr_1)$,$\tau_5(tr_1^2)$,$\tau_5(d_2^2\,p_1)$,$\tau_5(tr_1\,p_1\,si_1).$
\noindent
Put
$$
\begin{array}{ll}
c_1:=\tau_5(tr_2) &     [c_1]=[4,0,10],\\
c_2:=\tau_4(tr_3) &      [c_2]=[4,4,8],\\
c_3:=\tau_5(tr_3) &     [c_3]=[4,6,7],\\
s_1:=\tau_4(d_2\,tr_1) & [s_1]=[6,2,14],\\
v_1:=\tau_5(d_2\,p_2) & [v_1]=[8,0,20],\\
v_2:=\tau_5(d_2^2\,tr_1) &[si_1]=[8,2,19]\\
dv:=\tau_5(dv_2^2\,si_1)  & [dv]=[12,0,30]\\
vis:=\tau_5(p_1\,p_2\,si_1)  & [vis]=[18,0,45].
\end{array}
$$
\noindent
Thus 
$$ B_3=k[t,d_1,d_2,tr_1,tr_2,tr_3,c_1,c_2,c_3,p_1,p_2,s_1,si_1,v_1,v_2,dv,vis] $$
\noindent
The following polynomials are  acceptable for  $B_3$:
$$
\begin{array}{ll}
p_3:=\tau_2(c_3) & [p_3]=[5,7,9],\\
s_2:=\tau_1(p_1) & [s_2]=[6,4,13], \\
si_2:=\tau_1(s_1) & [si_2]=[7,5,15] ,\\
dev:=\tau_2(v_2) & [dev]=[9,3,21],\\
od:=\tau_5(s_1\,d_2^2) & [od]= [11,1,27],\\
trn:=\tau_5(v_2\,d_2^2) &  [trn]=[13,1,32] .
\end{array}
$$
\noindent
As above we can show 
$$B_4=k[t,d_{1-2},tr_{1-3},c_{1-3},p_{1-3},s_{1-2},si_{1-2},v_{1-2},dev,od,dv,trn,vis],$$
\noindent
Similarly we obtain 
$\tau(B_4) \subset B_4$ hence the indicated polynomial set  of 23 polynomials is a minimal generating set for  $k[X]^{\,d}.$

\noindent
{\bf{5.6 \rm {\bf n\,=\,6.}}}
\noindent
We have 
$B_2=k[t,d_1,d_2,d_3]$ where
$$
\begin{array}{ll}
d_1:=\tau_6(t) & [p_3]=[2,0,6],\\
d_2:=\tau_4(t) & [s_2]=[2,4,4], \\
d_3:=\tau_2(t) & [si_2]=[2,8,2] ,\\
\end{array}
$$
\noindent
Further $B_3=k[t,d_{1-3},tr_{1-4},p_{1-2}]$ where 
$$
\begin{array}{ll}
tr_1:=\tau_6(d_3) & [p_3]=[3,2,8],\\
tr_2:=\tau_4(d_3) & [s_2]=[3,6,6], \\
tr_3:=\tau_3(d_3) & [si_2]=[3,8,6] ,\\
tr_4:=\tau_1(d_3) & [si_2]=[3,12,3] ,\\
p_1:=\tau_6(d_2^{\,2}) & [p_1]=[5,2,14],\\
p_2:=\tau_5(d_2^{\,2}) & [p_1]=[5,4,13].\\
\end{array}
$$
\noindent
In the same way we obtain
$B_4=k[t,d_{1-3},tr_{1-4},c_{1-4},p_{1-2},s_{1-3},si_{1-2},vi,dev,de_{1-2},dvan]$ where 
$$
\begin{array}{ll}
c_1:=\tau_6(tr_2) & [c_1]=[4,0,12],\\
c_2:=\tau_5(tr_3) & [c_2]=[4,4,10], \\
c_3:=\tau_6(tr_4) & [c_3]=[4,6,9] ,\\
c_4:=\tau_4(tr_4) & [c_4]=[4,10,7] ,\\
s_1:=\tau_6(d_2\,tr_1) & [p_1]=[6,0,18],\\
s_2:=\tau_3(d_2\,tr_1) & [p_1]=[6,6,15],\\
s_3:=\tau_1(p_1) & [p_1]=[6,6,15],\\
si_1:=\tau_4(tr_1^2) & [si_1]=[7,2,20],\\
si_2:=\tau_3(tr_1^2) & [si_2]=[7,4,19],\\
vi:=\tau_6(d_2\,p_2) & [vi]=[8,2,23],\\
dev:=\tau_4(tr_1 p_2) & [dev]=[9,4,25],\\
de_1:=\tau_6(tr_1^3) &  [de_1]=[10,0,30],\\
de_2:=\tau_5(tr_1^3) &  [de_2]=[10,2,29],\\
dvan:=\tau_5(tr_1^2 p_1) & [dvan]=[12,2,35].
\end{array}
$$
\noindent
Acceptable polynomials for  $B_3$ are only these two polynomials :
$$
\begin{array}{ll}
p_3:=\tau_4(c_4) &  [p_3]=[5,8,11],\\
pt:=\tau_6(si_1\, si_2) & [pt]=[15,0,45].
\end{array}
$$
\noindent
The polynomial  $pt$ has  degree 15  order 0 and consist of 1370 terms. A weight of $pt$ is an odd number whereas the weights of all other generating polynomials of weight zero  are even numbers.  Therefore $pt$ is irreducible.
\noindent
We can show that 
$$B_4=k[t,d_{1-3},tr_{1-4},c_{1-4},p_{1-3},s_{1-3},si_{1-2},vi,dev,de_{1-2},dvan,pt].$$
\noindent
As above we may obtain  $\tau(B_4) \subset B_4 $  hence the indicated polynomial set  of 26 polynomials is a minimal generating set for  $k[X]^{\,d}.$

\noindent
{\bf Reference}
\begin{enumerate}

\item[1.] A. Nowicki. Polynomial derivation and their Ring of Constants. 
UMK, Torun,1994.

\item[2.]
 A. Nowicki. The fourteenth problem of Hilbert for polynomial 
derivations.Differential Galois theory. Proceedings of the 
workshop, Bedlewo, Poland, May 28-June 1, 2001. Warsaw: Polish 
Academy of Sciences, Institute of Mathematics, Banach Cent. Publ. 
58, 177-188, 2002. 

\item[3.]
 A. van den Essen. Polynomial automorphisms and the Jacobian 
conjecture. Progress in Mathematics (Boston, Mass.). 190. 
Basel.2000.

\item[4.] L.P. Bedratyuk. Symmetrical invariants of modular Lie algebras. Ph.D thesis (Russian), - Moscow state university,  1995.

\item[5.] A. Cerezo.  Tables des invariants alg$\acute{\mbox e}$briques et rationnels d'une matrice nilpotente de petite dimension. Pr$\acute{\mbox e}$publications Math$\acute{ \mbox e}$matiques, Universit$\acute {\mbox e}$ de Nice, 146, 1987.

\item[6.] R. Weitzenb$\ddot{\rm{o}}$ck   $\ddot{\rm{U}}$ber die Invarianten von linearen Gruppen. Acta Math.  58, 231-293, 1932.

\end{enumerate}

\end{document}